\documentclass[12pt]{article}
\usepackage{amsmath}
\usepackage{amssymb}
\usepackage{vatola}

\def\q{\quad}

\def\qtq#1{\q\t{#1}\q}
\def\mod#1{\ (\text{\rm mod}\ #1)}
\def\t{\text}
\def\f{\frac}
\def\e{\equiv}
\def\b{\binom}

\def\sls#1#2{(\f{#1}{#2})}
 
\def\Ls#1#2{\Big(\f{#1}{#2}\Big)}

\let \pro=\proclaim
\let \endpro=\endproclaim
\begin{document}

\par\q\par\q
 \centerline {\bf
Congruences for Catalan-Larcombe-French numbers}
$$\q$$
\centerline{Xiao-Juan Ji$^1$ and Zhi-Hong Sun$^2$}
\par\q\newline
\centerline{$\ ^1$School of Mathematical Sciences, Soochow
University,}
 \centerline{Suzhou, Jiangsu 215006,
 P.R. China}
\centerline{Email: Xiaojuanji2014@163.com}
\par\q\newline
\centerline{$\ ^2$School of Mathematical Sciences, Huaiyin Normal
University,} \centerline{Huaian, Jiangsu 223001, P.R. China}
\centerline{Email: zhihongsun@yahoo.com} \centerline{Homepage:
http://www.hytc.edu.cn/xsjl/szh}

 \abstract{Let $\{P_n\}$ be the Catalan-Larcombe-French numbers
  given by $P_0=1,\ P_1=8$ and
$n^2P_n=8(3n^2-3n+1)P_{n-1}-128(n-1)^2P_{n-2}$ $(n\ge 2)$, and let
$S_n=P_n/2^n$. In this paper we deduce congruences for
$S_{mp^r}\pmod{p^{r+2}}$, $S_{mp^r-1}\pmod{p^r}$ and
$S_{mp^r+1}\pmod{p^{2r}}$, where $p$ is an odd prime and $m,r$ are
positive integers. We also prove that $S_{(p^2-1)/2}\equiv 0\pmod
{p^2}$ for any prime $p\equiv 5,7\pmod 8$, and show that $\{S_m\}$
is log-convex.

 \par\q
 }
 \endabstract

\section*{1. Introduction}

\par Let $\{P_n\}$ be the sequence given by
$$P_0=1,\ P_1=8\qtq{and}
(n+1)^2P_{n+1}=8(3n^2+3n+1)P_n-128n^2P_{n-1}\ (n\ge 1).\tag 1.1$$
The numbers $P_n$ are called Catalan-Larcombe-French numbers since
Catalan first defined $P_n$ in [C], and in [LF1] Larcombe and French
proved that
$$P_n=2^n\sum_{k=0}^{[n/2]}(-4)^k\b{2n-2k}{n-k}^2\b{n-k}k
=\sum_{k=0}^n\f{\b{2k}k^2\b{2n-2k}{n-k}^2}{\b nk},\tag 1.2$$ where
$[x]$ is the greatest integer not exceeding $x$.
 The numbers $P_n$
occur in the theory of elliptic integrals, and are related  to the
arithmetic-geometric-mean. See [LF1] and A053175 in Sloane's
database ``The On-Line Encyclopedia of Integer Sequences".
 \par Let $\{f_n\}$ be the Franel numbers given by
 $f_n=\sum_{k=0}^n\b nk^3\q (n=0,1,2,\ldots)$, and let $\{S_n\}$
be defined by
$$S_0=1,\ S_1=4\qtq{and}
(n+1)^2S_{n+1}=4(3n^2+3n+1)S_n-32n^2S_{n-1}\ (n\ge 1).\tag 1.3$$
Comparing (1.3) with (1.1), we see that
$$S_n=\f{P_n}{2^n}.\tag 1.4$$
Zagier noted that
$$S_n=\sum_{k=0}^{[n/2]}\binom{2k}k^2\binom{n}{2k}4^{n-2k}.
\tag 1.5$$
 In this paper we investigate the properties of $S_n$
instead of $P_n$ since $S_n$ is an Ap\'epy-like sequence. As
observed by V. Jovovic in 2003 (see [LF2]),
$$ S_n=\sum_{k=0}^n\binom nk\binom
{2k}k\binom{2n-2k}{n-k}\q (n=0,1,2,\ldots). \tag 1.6$$ Recently Z.W.
Sun stated that
$$S_n=\sum_{k=0}^n\b{2k}k^2\b k{n-k}(-4)^{n-k}=\f 1{(-2)^n} \sum_{k=0}^n\b{2k}k\b{2n-2k}{n-k}
\b k{n-k}(-4)^k.\tag 1.7$$
 The first few values of $S_n$ are shown
below:
$$\align &S_0=1,\ S_1=4,\ S_2=20,\ S_3=112,\ S_4=676,\ S_5=4304,
\\&S_6=28496,\ S_7=194240,\ S_8=1353508,\ S_9=9593104.\endalign$$
\par Let $p$ be an odd prime. In [JLF], Jarvis, Larcombe and French
proved that if $n=a_rp^r+\cdots+a_1p+a_0$ with $a_0,a_1,
\ldots,a_r\in\{0,1,\ldots,p-1\}$, then
$$P_n\e P_{a_r}\cdots P_{a_1}P_{a_0}\mod p.\tag 1.8$$
 In [JV] Jarvis and Verrill showed that
$$P_n\e (-1)^{\f{p-1}2}128^nP_{p-1-n}\mod p\qtq{for}
n=0,1,\ldots,p-1\tag 1.9$$ and
$$P_{mp^r}\e P_{mp^{r-1}}\mod{p^r}\qtq{for}m,r\in\Bbb Z^+,\tag
1.10$$ where $\Bbb Z^+$ is the set of positive integers. In [OS]
Osburn and Sahu stated that
$$S_{mp^r}\e S_{mp^{r-1}}\mod{p^{2r}}\qtq{for}m,r\in\Bbb Z^+.\tag
1.11$$  But they did not give the details for the proof. In this
paper we will prove (1.11) in an elementary and natural manner. Let
$\varphi(n)$ be the Euler's totient function. Since $P_n=2^nS_n$,
from (1.11) we deduce a congruence for
$P_{mp^r}-2^{m\varphi(p^r)}P_{mp^{r-1}} \mod{p^{2r}}$, which
improves (1.10). Thus (1.11) is a vast generalization of (1.10).

\par In [S3] the second author established
 some identities involving $S_n$. For
example,
$$\sum_{k=0}^n\b nk(-1)^k\f{S_k}{8^k}=\f{S_n}{8^n}
\qtq{and}
\sum_{k=0}^{2n}\b{2n}k\b{2n+k}k(-8)^{2n-k}S_k=(-1)^n\b{2n}n^3.\tag
1.12$$
  For a prime
$p$ let $\Bbb Z_p$ denote the set of those rational numbers whose
denominator is not divisible by $p$.  Let $p$ be an odd prime,
$n\in\Bbb Z_p$ and $n\not\e 0,-16\mod p$. In [S3] the second author
proved that
$$\sum_{k=0}^{p-1}\b{2k}k\f{S_k}{(n+16)^k}\e
\Ls {n(n+16)}p\sum_{k=0}^{p-1} \f{\b{2k}k^2\b{4k}{2k}}{n^{2k}}\mod
p,\tag 1.13$$ where $\sls ap$ is the Legendre symbol.
\par Let
$r\in\Bbb Z^+$ and $ p$ be a prime with $p\e 5,7\pmod 8$. In [S3]
the second author conjectured that $$S_{\f{p^r-1}2}\e0\pmod
{p^r}\qtq{and}f_{\f{p^r-1}2}\e0\pmod {p^r}.\tag 1.14$$ In this paper
we prove (1.14) in the case $r=2$.

\par Let
 $\{E_n\}$ be the Euler numbers given by
$$E_{2n-1}=0, \ E_0=1\qtq{and} \sum_{k=0}^n\b{2n}{2k}E_{2k}=0\ (n\ge
1).$$ Suppose that $p$ is an odd prime and $n\in\Bbb Z^+$. In this
paper we determine $S_{mp^r}\mod {p^{r+2}}$ by showing that
$$S_{np}-S_n
\e\cases 8n^2S_{n-1}(-1)^{\f{p-1}2}p^2E_{p-3}\mod{p^3}&\t{if $p>3$
and $p\nmid n$,}
\\9(n-1)S_n\mod {p^3}&\t{if $p=3$
and $3\nmid n$,}
\\0 \mod {p^{3+\t{\rm
ord}_pn}}&\t{if $p\mid n$,}\endcases\tag 1.15$$ where
$\text{ord}_pn$ is the unique nonnegative integer $\alpha$ such that
$p^{\alpha}\mid n$ and $p^{\alpha+1}\nmid n$. We also show that
$$S_{mp^r+1}\e 4(mp^r+1)S_{mp^{r-1}}\mod{p^{2r}}\ \t{and}\
S_{mp^r-1}\e (-1)^{\f{p-1}2}S_{mp^{r-1}-1}\mod{p^{r}}.\tag 1.16$$ In
Section 4 we prove the second author's conjecture (see [S3])
$$S_m^2<S_{m+1}S_{m-1}<\big(1+\f 1{m(m-1)}\big)S_m^2\qtq{for}m=2,3,\ldots.$$
\section*{2. Basic lemmas}
  \pro{Lemma 2.1 (Lucas theorem [M])} Let $p$ be an odd prime.
 Suppose $a=a_rp^r+\cdots+a_1p+a_0$ and $b=b_rp^r+\cdots+b_1p+b_0$,
 where $a_r,\ldots,a_0,b_r,\ldots,b_0\in\{0,1,\ldots,p-1\}$. Then
 $$\b ab\e \b{a_r}{b_r}\cdots \b{a_0}{b_0}\mod p.$$
 \endpro
 \pro{Lemma 2.2} Let $p$ be an odd prime and $a,b\in\Bbb Z^+$.
  Suppose $a_0,b_0\in\{0,1,\ldots,p-1\}$. Then
  $$\b{ap+a_0}{bp+b_0}\e \b ab\b{a_0}{b_0}\mod p.$$
  \endpro
  Proof. Assume $a=a_rp^r+\cdots+a_2p+a_1$ and $b=b_rp^r+
  \cdots+b_2p+b_1$,
  where $a_r,\ldots,a_1,b_r,\ldots,b_1\in\{0,1,\ldots,p-1\}$. By
  Lucas theorem,
  $$\b{ap+a_0}{bp+b_0}\e \b{a_r}{b_r}\cdots\b{a_1}{b_1}\b{a_0}{b_0}
  \e \b ab\b{a_0}{b_0}\mod p.$$
  This is the result.
  \pro{Lemma 2.3 (Kazandzidis' congruence [M])} Let $p\ge 5$ be a prime and $m,n\in\Bbb Z^+$.
  Then
  $$\b{mp}{np}\e \b mn\mod{p^3}.$$
  \endpro

 \pro{Lemma 2.4 ([Su, Lemma 2.1])} Let $p$ be an odd prime and $k\in\{1,2,\ldots,p-1\}$.
 Then
 $$\b{2k}k\b{2(p-k)}{p-k}\e \cases -\f{2p}k\mod{p^2}&\t{if $k<\f
 p2$,}
 \\\f{2p}k\mod{p^2}&\t{if $k>\f p2$.}
 \endcases$$
 \endpro
\par Let $\{B_n\}$ be the Bernoulli numbers defined by
$B_0=1$ and $\sum_{k=0}^{n-1}\b nkB_k=0$ $(n\ge 2)$. It is known
that $B_{2k+1}=0$ for $k\in\Bbb Z^+$. For $m,n\in\Bbb Z^+$ it is
well known that
$$\sum_{k=0}^{n-1}k^m=\f
1{m+1}\sum_{k=1}^{m+1}\b{m+1}kB_kn^{m+1-k}.\tag 2.1$$ By the
Staudt-Clausen theorem, $B_{2k}\in \Bbb Z_p$ for $2k\not\e 0\mod
{p-1}$, and $pB_{2k}\in \Bbb Z_p$ for $2k\e 0\mod {p-1}$. See [MOS].
\par Let
$\{E_n(x)\}$ be the Euler polynomials given by
$$E_n(x)=\f 1{2^n}\sum_{k=0}^n\b nk(2x-1)^{n-k}E_k.$$
Then $E_n=2^nE_n(\f 12)$. It is known that (see [MOS])
$$\sum_{k=0}^{n-1}(-1)^kk^m=\f{E_m(0)-(-1)^nE_m(n)}2.$$

 \pro{Lemma 2.5 ([S2, Lemma 2.2])} Let $p$ be an odd
prime, $a\in\Bbb Z_p$, $a\not\e 0\mod p$ and
$k\in\{1,2,\ldots,p-2\}$. Then
$$\sum_{r=1}^{\langle a\rangle_p}\f{(-1)^r}{r^k}
\e -\f{(2^{p-k}-1)B_{p-k}}{p-k}+\f 12(-1)^{\langle
a\rangle_p+k}E_{p-1-k}(-a) \mod p,$$ where $\langle
a\rangle_p\in\{0,1,\ldots,p-1\}$ is given by $a\e \langle
a\rangle_p\mod p$.
\endpro

\pro{Lemma 2.6} Let $p$ be an odd prime, $k,m\in\Bbb Z^+$ and $\f
p2<k<p$. Then
$$\b{2mp+2k}{mp+k}\e (2m+1)\b{2m}m\b{2k}k\mod{p^2}.$$
\endpro
Proof. Clearly
$$\align &\f{(2m+1)p((2m+1)p-1)\cdots
((m+1)p+1) (m+1)p}{(mp)!}
\\&=\f{((2m+1)p)(2mp)\cdots ((m+1)p)}{p\cdot (2p)\cdots (mp)}
\cdot\f{\prod_{r=m+1}^{2m}(rp+1)\cdots(rp+p-1)}{\prod
_{r=0}^{m-1}(rp+1)\cdots(rp+p-1)}
\\&=p(2m+1)\b{2m}m\cdot\f{\prod_{r=m+1}^{2m}(rp+1)\cdots(rp+p-1)}{\prod
_{r=0}^{m-1}(rp+1)\cdots(rp+p-1)} \\&\e
p(2m+1)\b{2m}m\f{(p-1)!^m}{(p-1)!^m} =p(2m+1)\b{2m}m\mod{p^2}.
\endalign$$
Thus
$$\align
&\b{2mp+2k}{mp+k}
\\&=\f{(2m+1)p((2m+1)p-1)\cdots
((m+1)p+1) (m+1)p}{(mp)!}
\\&\q\times\f{(2mp+2k)\cdots(2mp+p+1)((m+1)p-1)\cdots ((m+1)p-(p-1-k))}
{(mp+1)\cdots(mp+k)}
\\&\e p(2m+1)\b{2m}m\f{(2k)(2k-1)\cdots(p+1)(p-1)(p-2)\cdots(k+1)}
{k!}\\&=(2m+1)\b{2m}m\b{2k}k\mod{p^2}.
\endalign$$
This proves the lemma.
 \pro{Lemma 2.7} For any positive integer $n$
we have
$$S_n=2\sum_{k=1}^n\b{n-1}{k-1}\b{2k}k\b{2n-2k}{n-k}.$$
\endpro
Proof. Since
$$\align \sum_{k=0}^n(2k-n)\b nk\b{2k}k\b{2n-2k}{n-k}
&=\sum_{k=0}^n (2(n-k)-n)\b nk\b{2k}k\b{2n-2k}{n-k}
\\&=-\sum_{k=0}^n(2k-n)\b nk\b{2k}k\b{2n-2k}{n-k},\endalign$$
we see that $$\sum_{k=0}^n(2k-n)\b nk\b{2k}k\b{2n-2k}{n-k}=0\tag
2.2$$ and so
$$nS_n=\sum_{k=0}^n2k\b nk\b{2k}k\b{2n-2k}{n-k}=2n\sum_{k=1}^n
\b{n-1}{k-1}\b{2k}k\b{2n-2k}{n-k}.$$ This yields the result.

\pro{Lemma 2.8} Let $m\in \Bbb Z$ and $k,n,p\in \Bbb Z^+$. Then
$$\b {m{p^r}-1}{k}=(-1)^{k-[\f kp]}\b {m{p^{r-1}}-1}{[k/p]}
\prod\Sb i=1\\p\nmid i\endSb^k\Big(1-\f{mp^r}i\Big)$$ and so
$$\b {m{p^r}}{np}=(-1)^{n(p-1)}\b {m{p^{r-1}}}n
\prod\Sb {i=1}\\p \nmid i\endSb^{np-1}\Big(1-\f {mp^r}i\Big).$$
\endpro
Proof.
 Clearly,
$$\aligned\b{mp^r-1}k&=\prod_{i=1}^{k}\f{mp^r-i}i=\prod\Sb i=1
\\p\nmid i\endSb^k\f{mp^r-i}i
\prod_{i=1}^{[k/p]}\f{mp^r-pi}{pi}
\\&=\prod\Sb i=1\\p\nmid i\endSb^k\f{mp^r-i}i
\prod_{i=1}^{[k/p]}\f{mp^{r-1}-i}i
 \\&=(-1)^{k-[\f kp]}\prod\Sb
i=1\\p\nmid
 i\endSb^k\Big(1-\f{mp^r}i\Big)
\cdot \b {m{p^{r-1}}-1}{[k/p]}.\endaligned$$
 Taking $k=np-1$ in the above we see that
$$\b{mp^r}{np}=\f{mp^r}{np}\b{mp^r-1}{np-1}
=\f{mp^{r-1}}n\b{mp^{r-1}-1}{n-1}\cdot (-1)^{np-1-(n-1)}\prod\Sb
{i=1}\\p \nmid i\endSb^{np-1}\Big(1-\f {mp^r}i\Big).$$ This yields
the remaining part.

\pro{Lemma 2.9} Let $m,n\in\Bbb Z^+$.Then
$$\b{3m}{3n}\e \b mn(1+9mn^2-9m^2n)\mod{27}.$$ \endpro
 Proof. By Lemma 2.8,
 $$\aligned\b{3m}{3n}&=\b mn\prod\Sb k=1\\3\nmid k\endSb ^{3n-1}
 \big(1-\f{3m}k\big)
\\&\e\b mn\Big( 1-3m\sum\Sb k=1\\3\nmid k\endSb ^{3n-1}\f 1k
+9m^2\sum\Sb 1\le i<j\le {3n-1}\\3\nmid {ij}\endSb\f 1{ij}\Big)
\pmod{27}.\endaligned$$
 Clearly
 $$\aligned 2\sum\Sb 1\le i<j\le {3n-1}\\3\nmid {ij}\endSb\f 1{ij}
&=\Big(\sum\Sb i=1\\3\nmid i\endSb ^{3n-1}\f 1i\Big)^2
 -\sum\Sb i=1\\3\nmid {i}\endSb ^{3n-1}\f 1{i^2}
 \e \Big(\sum _{i=1}^{3n-1}i\Big)^2
 -\sum\Sb i=1\\3\nmid i\endSb ^{3n-1}1
\\&=\Big( \f{(3n-1)3n}2  \Big)^2-(3n-1-(n-1))
\e -2n \pmod 3.\endaligned$$
 Thus,
$$
\b{3m}{3n}\e\b mn\Big(1- 3m\sum \Sb k=1\\3\nmid {k}\endSb ^{3n-1}\f
1k-9m^2n \Big)\pmod {27}.$$
 By Euler's Theorem, for $k\not\e 0 \pmod 3$
  we have $k^6=k^{\varphi (9)}\e 1\pmod 9.$
 Thus
 $$\aligned\sum \Sb k=1\\3\nmid k\endSb^{3n-1}\f 1k
 &\e\sum_{k=1}^{3n-1}k^5=\f{B_6(3n)-B_6}6=\f 16\sum_{k=1}^6\b
6k(3n)^k{B_{6-k}}
\\&\e\f 16\b 62(3n)^2{B_4}=-\f 34 n^2\e-3n^2\pmod 9.\endaligned$$
 Hence
 $$\aligned\b {3m}{3n}&\e\b mn(1+9mn^2-9m^2n)
 \\&\e\cases\b mn \pmod {27}  &\text{if $3\mid mn(m-n) $ ,}
 \\
  \b mn (1+9n)\pmod {27}  &\text{if $3\nmid mn$ and $3\mid {m+n}
 $}.\endcases\endaligned $$
\pro{Lemma 2.10} Let $p$ be an odd prime, $r,m\in\Bbb Z^+$ and
$s\in\{0,1,\ldots,mp^{r-1}-1\}$. Then
$$\b{mp^r}{sp}\e\b{mp^{r-1}}s\pmod{p^{2r}}.$$
\endpro
 Proof. By Lemma 2.8, $$\b{mp^r}{sp}=
\b{mp^{r-1}}s\prod\limits\Sb {i=1}\\p\nmid i\endSb ^{sp-1}
\Big(1-\f{mp^r}i\Big) \e\b{mp^{r-1}}s\Big(1-mp^r\sum\limits\Sb
{i=1}\\p\nmid i\endSb ^{sp-1}\f 1i\Big)\pmod{p^{2r}}.$$ Set
$l=\t{ord}_ps$ and $s=p^ls_0.$
 Since $B_1=-\f 12$, $B_{2s+1}=0\ (s\ge 1)$, $pB_k\in\Bbb Z_p$ and
 $\varphi(p^{l+1})\ge l+2$  we see that
$$\aligned \sum\Sb
i=1\\p\nmid i\endSb^{sp-1}\f 1i&\e \sum\Sb i=1\\p\nmid
i\endSb^{sp-1}i^{\varphi(p^{l+1})-1} \e \sum_{i=1}^{sp-1}
i^{\varphi(p^{l+1})-1}
\\&=\f 1{\varphi(p^{l+1})}
\sum_{j=1}^{\varphi(p^{l+1})}
\b{\varphi(p^{l+1})}j(s_0p^{l+1})^jB_{\varphi(p^{l+1})-j}
\\&=(s_0p^{l+1})^{\varphi(p^{l+1})-1}B_1+\f
1{\varphi(p^{l+1})}\sum_{j=1}^{\varphi(p^{l+1})/2}
\b{\varphi(p^{l+1})}{2j} (s_0p^{l+1})^{2j}B_{\varphi(p^{l+1})-2j}
\\&\e \f {s_0}{p-1}\sum_{j=1}^{\varphi(p^{l+1})/2}
\b{\varphi(p^{l+1})}{2j}(s_0p^{l+1})^{2j-1} \cdot
pB_{\varphi(p^{l+1})-2j}\e 0\mod{p^{l+1}}.\endaligned\tag 2.3$$
 If $l\ge r-1$, then $r+l+1\ge 2r$ and so
 $$\b{mp^r}{sp}\e\b{mp^{r-1}}s\Big(1-mp^r\sum\limits\Sb
 {i=1}\\p\nmid i\endSb
^{sp-1}\f 1i\Big)\e\b{mp^{r-1}}s\pmod {p^{2r}}.$$
 If $0\le l<r-1$, then $\b{mp^{r-1}}s
 =\f{mp^{r-1}}{s_0p^l}\b{mp^{r-1}-1}{s-1}
 \e 0\pmod{p^{r-1-l}}$ and so
$$\b{mp^r}{sp}\e\b{mp^{r-1}}s-mp^r\b{mp^{r-1}}s
\sum\limits\Sb {i=1}\\p\nmid i\endSb
 ^{sp-1}\f 1i\e\b{mp^{r-1}}s\pmod {p^{2r}}.$$
 This completes the proof.

 \pro{Lemma 2.11} Let $p$ be an odd prime, $r,m\in\Bbb Z^+$ and
$s\in\{1,2,\ldots,mp^{r-1}\}$. Then
$$\align&\b{mp^{r-1}}s\b{2sp}{sp}\b{2(mp^{r-1}-s)p}{(mp^{r-1}-s)p}
\\&\e\cases \b ms\b {2s}s\b {2(m-s)}{m-s}(1+9m)\pmod {p^{r+2}}
 &\text {if $r=1\ and\ p=3$,}
\\\b{mp^{r-1}}s\b{2s}s\b{2(mp^{r-1}-s)}{mp^{r-1}-s}\mod{p^{r+2}}
 &\text {if $r>1$ or $p>3$.}\endcases\endalign$$
\endpro
Proof. For $r=1$ the result follows from Lemmas 2.3 and 2.9.
 Now assume $r\ge
2$. If $p\nmid s$, then
$\b{mp^{r-1}}s=\f{mp^{r-1}}s\b{mp^{r-1}-1}{s-1}
\e 0\mod{p^{r-1}}$.
By Lemmas 2.3 and 2.9,
$$\b{2sp}{sp}\b{2(mp^{r-1}-s)p}{(mp^{r-1}-s)p}
\e \b{2s}s\b{2(mp^{r-1}-s)}{mp^{r-1}-s}\mod{p^3}.$$
Thus the result is true.
 Now assume that $p\mid s$, $l=\t{ord}_ps$
and $s=p^ls_0$. For $1\le l<r-1$, using Lemma 2.10 we see that
$$\b{2sp}{sp}\b{2(mp^{r-1}-s)p}{(mp^{r-1}-s)p}
\e \b{2s}s\b{2(mp^{r-1}-s)}{mp^{r-1}-s}\mod{p^{2l+2}}.$$ Since
$\b{mp^{r-1}}s=\f{mp^{r-1}}{p^ls_0}\b{mp^{r-1}-1}{s-1}\e
0\mod{p^{r-1-l}}$ and $r-1-l+2l+2=r+l+1\ge r+2$, the result is true
in this case.
 For $l\ge r$ we see that $p^r\mid sp$ and $p^r\mid
(mp^{r-1}-s)p$. Thus applying Lemma 2.10 we deduce that
$$\b{2sp}{sp}\b{2(mp^{r-1}-s)p}{(mp^{r-1}-s)p}
\e \b{2s}s\b{2(mp^{r-1}-s)}{mp^{r-1}-s}\mod{p^{2r}}.$$
 As $2r\ge r+2$, the result is also true.
  The proof is now complete.

 \pro{Lemma 2.12} Let $p$ be an odd prime, $m,r\in\Bbb Z^+ $ and
 $k\in \{0,\ldots,mp^r\}$.
 Then
 $$k\b{2k}k\b{2(mp^r-k)}{mp^r-k}\e 0\pmod {p^r}.$$
 \endpro
 Proof. Suppose $s=[\f kp]$ and $t=k-sp$. Then $t\in\{0,1,\ldots,
 p-1\}$. We first assume $p\nmid k$.
    That is, $t>0$.
      Let us consider the case $r=1$.
    By Lemma 2.2,
 for $1\le t<\f p2$,
  $$\align \b{2k}k\b{2(mp-k)}{mp-k}&=\b{2k}k
  \b{(2(m-s)-1)p+
  p-2t}
 {(m-s-1)p+p-t}
  \\&\e \b{2k}k\b{2(m-s)-1}
  {m-s-1}\b{p-2t}{p-t}=0
  \mod p,\endalign$$
  and for $t>\f p2$,
$$\align \b{2k}k\b{2(mp-k)}{mp-k}&=\b{(2s+1)p+2t-p}
{sp+t}\b{2(mp-k)}{mp-k}
    \\&\e \b{2s+1}{s}\b{2t-p}t\b{2(mp-k)}{mp-k}=0
  \mod p.\endalign$$ Thus the result is true for $r=1$.
\par Now assume $p\nmid k$ and $r\ge 2$. Suppose that for $n<r$ and
 $k\in\{1,2,\ldots,mp^n-1\}$ we have
$$\b{2k}k\b{2(mp^n-k)}{mp^n-k}\e 0\pmod
{p^n}.$$ When $p\mid s$,  by the inductive hypothesis we have
$$\aligned&\b{2s}s\b{2mp^{r-1}-2s-2}{mp^{r-1}-s-1}(2mp^{r-1}-2s-1)p
\\&=\f{s+1}{2(2s+1)}(2mp^{r-1}-2s-1)p
\b{2s+2}{s+1}\b{2mp^{r-1}-2s-2}{mp^{r-1}-s-1}
 \e 0\pmod {p^r}.\endaligned$$
When $p\nmid s$, by the inductive hypothesis  we obtain
$$\aligned &\b{2s}s\b{2mp^{r-1}-2s-2}{mp^{r-1}-s-1}(2mp^{r-1}-2s-1)p
\\&=\f{mp^{r-1}-s}2p \b{2s}{s}\b{2mp^{r-1}-2s}{mp^{r-1}-s}
 \e 0 \pmod{p^r}.\endaligned$$
\par Suppose $k\in\{1,2,\ldots,mp^r-1\}$. For $t< \f {p}2$, from the
above we see that
 $$\aligned&\b{2k}k\b{2(mp^r-k)}{mp^r-k}
 =\b{2sp+2t}{sp+t}\b{(2mp^{r-1}-2s-1)p+p-2t}{(mp^{r-1}-s-1)p+p-t}
\\&=\b{2s}s\b{2mp^{r-1}-2s-2}{mp^{r-1}-s-1}(2mp^{r-1}-2s-1)pQ\e 0
\pmod {p^r},
\endaligned$$
and for $t> \f {p}2$,
 $$\aligned\b{2k}k\b{2(mp^r-k)}{mp^r-k}
 &=\b{(2s+1)p+2t-p}{sp+t}
 \b{(2mp^{r-1}-2s-2)p+2p-2t}{(mp^{r-1}-s-1)p+p-t}
\\&=Q\b{2s}s\b{2mp^{r-1}-2s-2}{mp^{r-1}-s-1}(2s+1)p
\\&\e -Q
\b{2s}s\b{2mp^{r-1}-2s-2}{mp^{r-1}-s-1}(2mp^{r-1}-2s-1)p
\\&\e 0\pmod{p^r},\endaligned$$
where
$$Q=\f{\sum\limits\Sb {i=1}\\p\nmid i \endSb
 ^{2sp+2t}i\sum\limits\Sb {i=1}\\p\nmid i \endSb
 ^{2(mp^r-sp-t)}i}{\sum\limits\Sb {i=1}\\p\nmid i \endSb
 ^{sp+t}i^2\sum\limits\Sb {i=1}\\p\nmid i \endSb
 ^{mp^r-sp-t}i^2}\in\Bbb Z_p.$$
Hence the result is true for $n=r$. Summarizing the above we prove
the result in the case $p\nmid k$.
 \par Now we assume $p\mid k$. Set $l=\t{ord}_p{k}$ and $k=p^lk_0$.
 Then $k_0\in \{1,\ldots,mp^{r-l}-1 \}$ and $p\nmid k_0$.
 For $l\ge r$ obviously we have
$k\b{2k}k\b{2(mp^r-k)}{mp^r-k} \e 0 \pmod {p^r}.$ For $1\le l\le
r-1$, since $p\nmid k_0$, from the above we deduce that
 $$ k\b{2k}k\b{2(mp^r-k)}{mp^r-k}
=Wp^l\b{2k_0}{k_0}\b{2mp^{r-l}-2k_0}{mp^{r-l}-k_0} \e0 \pmod
{p^r},$$
 where $W\in\Bbb Z_p$. The proof is now complete.

\pro{Lemma 2.13} Let $p$ be an odd prime, $r,m\in\Bbb Z^+$ and
$s\in\{0,1,\ldots,mp^{r-1}-1\}$. Then
$$\b{2sp+p-1}{sp+\f{p-1}2}
\b{2(mp^{r-1}-s-1)p+p-1}{(mp^{r-1}-s-1)p+\f{p-1}2}
\e\b{2s}s\b{2(mp^{r-1}-s-1)}{mp^{r-1}-s-1}\mod{p^r}.$$
\endpro
 Proof. Using Lemma 2.8 and the identity
 $$\b{a-b}{c-d}\b bd=\f{\b ac\b cd\b {a-c}{b-d}}{\b ab}\tag 2.4$$
 we have
 $$\aligned &\b{2sp+p-1}{sp+\f{p-1}2}
 \b{2(mp^{r-1}-s-1)p+p-1}{(mp^{r-1}-s-1)p+(p-1)/2}
 \\&=\f{mp^{r-1}}{2(2s+1)}\cdot
 \f{\b{2mp^r}{mp^r}\b{mp^r-1}{sp+(p-1)/2}^2}{\b{2mp^r-1}{2sp+p}}
\\&=\f{mp^{r-1}}{2(2s+1)}\cdot
 \f{\b{2mp^{r-1}}{mp^{r-1}}
 \b{mp^{r-1}-1}s^2}{\b{2mp^{r-1}-1}{2s+1}}\cdot
\f{\prod\limits \Sb i=1 \\p\nmid i\endSb ^{mp^r}\f {2mp^r-i}i
(\prod\limits \Sb i=1 \\p\nmid i\endSb ^{sp+(p-1)/2} \f{mp^r-i}i)^2}
{\prod\limits \Sb i=1 \\p\nmid i\endSb
^{2sp+p}\f{2mp^r-i}i},\endaligned$$ and
$$\b{2s}s\b{2(mp^{r-1}-s-1)}{mp^{r-1}-s-1}
=\f{mp^{r-1}}{2(2s+1)}\cdot
 \f{\b{2mp^{r-1}}{mp^{r-1}}
 \b{mp^{r-1}-1}s^2}{\b{2mp^{r-1}-1}{2s+1}}.$$
Since
$$\aligned\f{\b{2sp+p-1}{sp+(p-1)/2}
\b{2(mp^{r-1}-s-1)p+p-1}{(mp^{r-1}-s-1)p+(p-1)/2}}
{\b{2s}s\b{2(mp^{r-1}-s-1)}{mp^{r-1}-s-1}} &=\f{\prod\limits \Sb i=1
\\p\nmid i\endSb ^{mp^r}\f {2mp^r-i}i
{\prod\limits \Sb i=1 \\p\nmid i\endSb ^{sp+(p-1)/2}(\f
{mp^r-i}i})^2} {\prod\limits \Sb i=1
\\p\nmid i\endSb ^{2sp+p}\f {2mp^r-i}i}
 \e \f{\prod\limits \Sb i=1 \\p\nmid i\endSb ^{mp^r}(-1)
\prod\limits \Sb i=1 \\p\nmid i\endSb
^{sp+(p-1)/2}{(-1)^2}}{\prod\limits \Sb i=1 \\p\nmid i\endSb
^{2sp+p}(-1)}\\&=1\pmod {p^r},\endaligned$$ we obtain the result.

\pro{Lemma 2.14} Let $p$ be an odd prime, $r,m\in\Bbb Z^+$ and
$s\in\{0,1,\ldots,mp^{r-1}-1\}$. Then
$$\b{mp^r}{sp}\b{2sp}{sp}\b{2mp^r-2sp}{mp^r-sp}
\e\b{mp^{r-1}}s\b{2s}{s}\b{2mp^{r-1}-2s}{mp^{r-1}-s}\pmod{p^{2r}}.$$
\endpro
 Proof. Set $l=\t{ord}_ps$ and $s=p^ls_0.$
  By (2.3), (2.4), Lemmas 2.8 and 2.10,
 we have
 $$\aligned\f{\b{2sp}{sp}\b{2mp^r-2sp}{mp^r-sp}}
 {\b{2s}{s}\b{2mp^{r-1}-2s}{mp^{r-1}-s}}
&=\f{\prod\limits \Sb i=1\\p\nmid i\endSb ^{mp^r} \f {2mp^r-i}i
{\prod\limits \Sb i=1 \\p\nmid i\endSb ^{sp}
(\f{mp^r-i}i})^2}{\prod\limits \Sb i=1
\\p\nmid i\endSb ^{2sp}\f {2mp^r-i}i}
\e\f{\Big(-2mp^r\sum\limits \Sb i=1\\p\nmid i\endSb ^{mp^r}
 \f1i+1\Big)\Big(-2mp^r\sum\limits \Sb i=1
\\p\nmid i\endSb ^{sp}\f 1i+1\Big)}
{-2mp^r\sum\limits \Sb i=1
\\p\nmid i\endSb ^{2sp}\f 1i+1}
 \\& \e\f{-2mp^r\sum\limits \Sb i=1\\p\nmid i\endSb ^{mp^r}
 \f1i-2mp^r\sum\limits \Sb i=1
\\p\nmid i\endSb ^{sp}\f 1i+1}
{-2mp^r\sum\limits \Sb i=1
\\p\nmid i\endSb ^{2sp}\f 1i+1}\e 1 \pmod{p^{r+\min\{l+1,r\}}}.
\endaligned$$
  Thus,
 $$\b{2sp}{sp}\b{2mp^r-2sp}{mp^r-sp}
 \e\b{2s}{s}\b{2mp^{r-1}-2s}{mp^{r-1}-s}\pmod{p^{r+\min\{l+1,r\}}}.$$
\\If $l\ge r-1$, then $r+l+1\ge 2r$ and so
$$\b{mp^r}{sp}\b{2sp}{sp}\b{2mp^r-2sp}{mp^r-sp}
\e\b{mp^{r-1}}s\b{2s}{s}\b{2mp^{r-1}-2s}{mp^{r-1}-s}\pmod{p^{2r}}.$$
If $1\le l<r-1$, then $\b{mp^{r-1}}s\e0\pmod{p^{r-1-l}}$ and so
$$\b{mp^r}{sp}\b{2sp}{sp}\b{2mp^r-2sp}{mp^r-sp}
\e\b{mp^{r-1}}s\b{2s}{s}\b{2mp^{r-1}-2s}{mp^{r-1}-s}\pmod{p^{2r}}.$$
Now the proof is complete.

 \pro{Lemma 2.15 ([B])} Let $p$ be an odd prime.
  Suppose $n=n_1p+n_0$ and $k=k_1p+k_0$ with $k_1,n_1\in\Bbb Z^+$ and
  $k_0,n_0\in\{0,1,\ldots,p-1\}$. Then
  $$\b nk\e\b {n_1}{k_1}\Big[(1+n_1)\b{n_0}{k_0}-
  (n_1+k_1)\b{n_0-p}{k_0}-k_1\b{n_0-p}{k_0+p}\Big] \pmod {p^2}.$$
\endpro
\pro{Lemma 2.16} Let $p$ be an odd prime. Then
$$\sum_{t=0}^{(p-1)/2}(-1)^t
\Big[\b{\f {p-1}2+t}t-\b{p+\f{p-1}2+t}{p+t}\Big]
  \b{-\f{1}2}t^2\e 0 \pmod {p^2}.$$
\endpro
 Proof. From Lemma 2.2 we have
 $$\align\b{\f{p-1}2-p}{p+t}
 &=\b{-\f{p+1}2}{p+t}=(-1)^{t+1}\b{p+\f{p+1}2+t-1}{p+t}
 \e(-1)^{t+1}\b{\f{p+1}2+t-1}{t}
 \\&\e(-1)^{t+1}\b{\f{1}2+t-1}{t}
 =-\b{-\f{1}2}t
 \e-\b{\f{p-1}2-p}t
 \pmod p\endalign$$
 and so $\b{\f{p-1}2-p}t+\b{\f{p-1}2-p}{p+t}
 =(-1)^t\Big(\b{\f{p-1}2+t}t-\b{p+\f{p-1}2+t}{p+t}\Big)\e 0\pmod p.$
\par We first assume $p\e1\pmod 4$. Applying Lemma 2.15 we get
$$\align &\b{\f{3(p-1)}4}{\f{p-1}4}-\b{p+\f{3(p-1)}4}{p+\f{p-1}4}
=\b{\f{3(p-1)}4}{\f{p-1}2}-\b{p+\f{3(p-1)}4}{\f{p-1}2}
\\&\e\b{\f{3(p-1)}4}{\f{p-1}2}-\Big(2\b{\f{3(p-1)}4}{\f{p-1}2}
-\b{\f{3(p-1)}4-p}{\f{p-1}2}\Big)
 \\&=-\b{\f{3(p-1)}4}{\f{p-1}2}+(-1)^{\f{p-1}2}\b{\f{3(p-1)}4}{\f{p-1}2}
 =0\pmod {p^2},\endalign$$ and
 $$\aligned&\b{\f {p-1}2+t}t-\b{p+\f{p-1}2+t}{p+t}
 +\b{p-1-t}{\f{p-1}2-t}-\b{p+p-1-t}{p+\f{p-1}2-t}
 \\&=\b{\f {p-1}2+t}{\f {p-1}2}-\b{p+\f{p-1}2+t}{\f{p-1}2}
+\b{p-1-t}{\f{p-1}2}-\b{p+p-1-t}{\f{p-1}2}
 \\&\e\b{\f {p-1}2+t}{\f{p-1}2}
 -\Big(2\b{\f {p-1}2+t}{\f {p-1}2}-\b{\f{p-1}2-p+t}{\f{p-1}2} \Big)
 +\b{p-1-t}{\f{p-1}2}
 -\Big(2\b{p-1-t}{\f{p-1}2}-\b{-1-t}{\f {p-1}2}\Big)
\\&=-\b{\f {p-1}2+t}{\f {p-1}2}+(-1)^{\f{p-1}2}\b{p-1-t}{\f {p-1}2}
 -\b{p-1-t}{\f {p-1}2}+(-1)^{\f{p-1}2}\b{\f {p-1}2+t}{\f {p-1}2}
\\&=0\pmod {p^2}.\endaligned $$
Also,
$$\aligned
&(-1)^t\b{-\f{1}2}t^2-(-1)^{\f{p-1}2-t}\b{-\f{1}2}{\f{p-1}2-t}^2
 \e(-1)^t\b{-\f{1}2}t^2-(-1)^{t}\b{\f{p-1}2}{\f{p-1}2-t}^2
\\&\e(-1)^t\b{-\f{1}2}t^2-(-1)^t\b{-\f{1}2}{t}^2=0\pmod
p.\endaligned$$ Hence
 $$\align
 &\sum_{t=0}^{(p-1)/2}(-1)^t
 \Big[\b{\f {p-1}2+t}t-\b{p+\f{p-1}2+t}{p+t}\Big]
  \b{-\f{1}2}t^2
  \\&=\sum_{t=0}^{(p-5)/4}(-1)^t\b{-\f{1}2}t^2
  \Big[\b{\f{p-1}2+t}t-\b{p+\f{p-1}2+t}{p+t}\Big]
  \\&\q+
 (-1)^{\f{p-1}4}\b{-\f {1}2}{\f {p-1}4}^2
 \Big[\b{\f{3(p-1)}4}{\f{p-1}4}-\b{p+\f{3(p-1)}4}{p+\f{p-1}4}\Big]
 \\&\q+\sum_{t=0}^{(p-5)/4}(-1)^{\f{p-1}2-t}\b{-\f{1}2}{\f{p-1}2-t}^2
 \Big[\b{p-1-t}{\f{p-1}2-t}-\b{p+p-1-t}{p+\f{p-1}2-t}\Big]
\\&\e\sum_{t=0}^{(p-5)/4}(-1)^t\b{-\f{1}2}t^2
\Big[\b{\f{p-1}2+t}t-\b{p+\f{p-1}2+t}{p+t} \Big]
 \\&\q+(-1)^{\f{p-1}4}\b{\f {3(p-1)}4}
 {\f{p-1}4}
 \Big[\b{\f{3(p-1)}4}{\f{p-1}4}-\b{p+\f{3(p-1)}4}{p+\f{p-1}4}
\Big]
 \\&\q-\sum_{t=0}^{(p-5)/4}(-1)^{\f{p-1}2-t}\b{-\f{1}2}{\f{p-1}2-t}^2
 \Big[\b{\f{p-1}2+t}t-\b{p+\f{p-1}2+t}{p+t}\Big]
 \\&\e\sum_{t=0}^{(p-5)/4}\Big[(-1)^t\b{-\f{1}2}t^2
 -(-1)^{\f{p-1}2-t}\b{-\f{1}2}{\f{p-1}2-t}^2\Big]
 \Big[\b{\f{p-1}2+t}t-\b{p+\f{p-1}2+t}{p+t}\Big]
\\&\q+ (-1)^{\f{p-1}4}\b{\f {3(p-1)}4}{\f{p-1}4}
\Big[\b{\f{3(p-1)}4}{\f{p-1}4}
-\b{p+\f{3(p-1)}4}{p+\f{p-1}4}\Big]
  \e 0 \pmod {p^2}.\endalign$$
Thus the result is true for $p \e1\pmod 4$.
\par Now we assume $p \e3\pmod 4$. By Lemma 2.15,
$$\b{\f {p-1}2+t}t-\b{p+\f{p-1}2+t}{p+t}
\e
-\Big(\b{\f{p-1}2+t}{\f{p-1}2}+\b{p-1-t}{\f{p-1}2}\Big)\pmod{p^2}.$$
As $$\b{\f{p-1}2+t}{\f{p-1}2}+\b{p-1-t}{\f{p-1}2}
\e\b{\f{p-1}2+t}{\f{p-1}2}+\b{-1-t}{\f{p-1}2}
=\b{\f{p-1}2+t}{\f{p-1}2}+(-1)^{\f{p-1}2}\b{t+\f{p-1}2}{\f{p-1}2}
=0\pmod p$$ and $$\align&(-1)^t\b{-\f{1}2}t^2
 +(-1)^{\f{p-1}2-t}\b{-\f{1}2}{\f{p-1}2-t}^2
 \e(-1)^t\b{-\f{1}2}t^2-(-1)^{t}\b{\f{p-1}2}{\f{p-1}2-t}^2
\\&\e(-1)^t\b{-\f{1}2}t^2-(-1)^t\b{-\f{1}2}{t}^2=0\pmod p,
\endalign$$
 we obtain
$$\align&\sum_{t=0}^{(p-1)/2}
(-1)^t\Big[\b{\f {p-1}2+t}t-\b{p+\f{p-1}2+t}{p+t}\Big]
\b{-\f{1}2}t^2
\\&\e-\sum_{t=0}^{(p-3)/4}
 \Big[\b{\f{p-1}2+t}{\f{p-1}2}+\b{p-1-t}{\f{p-1}2}\Big]
  \Big[(-1)^t\b{-\f{1}2}t^2+(-1)^{\f{p-1}2-t}
  \b{-\f{1}2}{\f{p-1}2-t}^2\Big]
  \\&\e 0 \pmod {p^2}.
  \endalign$$
Hence the result is also true in this case. The proof is now
complete.

\section*{3. Congruences for
$S_{mp^r}\pmod{p^{r+2}}$, $S_{mp^r-1}\mod{p^r}$ and
$S_{mp^r+1}\mod{p^{2r}}$}

 \pro{Theorem 3.1} Let $p$ be an odd
prime and $n\in\Bbb Z^+$.  Then
$$S_{np}-S_n
\e\cases 8n^2S_{n-1}(-1)^{\f{p-1}2}p^2E_{p-3}\mod{p^3}&\t{if $p>3$
and $p\nmid n$,}
\\9(n-1)S_n\mod {p^3}&\t{if $p=3$
and $3\nmid n$,}
\\0 \mod {p^{3+\t{\rm
ord}_pn}}&\t{if $p\mid n$.}\endcases$$
\endpro
Proof. Set $r=\t{\rm ord}_p(np)$. Then
$$\align S_{np}&=\sum_{k=0}^{np}\b {np}k\b{2k}k\b{2(np-k)}{np-k}
\\&=\sum_{s=0}^n\b{np}{sp}\b{2sp}{sp}\b{2(n-s)p}{(n-s)p}
+\sum_{t=1}^{p-1}\sum_{s=0}^{n-1}\b{np}{sp+t}\b{2(sp+t)}{sp+t}
\b{2(np-sp-t)}{np-sp-t}.
\endalign$$
If $p>3$ or if $p=3$ and $3\mid n$, using Lemmas 2.10 and 2.11 we
see that $\b{np}{sp}\e \b ns\mod{p^{r+2}}$ and $$\b
ns\b{2sp}{sp}\b{2(n-s)p}{(n-s)p} \e \b ns\b{2s}s\b{2(n-s)}{n-s}
\mod{p^{r+2}}.$$ Thus,
$$\align &\sum_{s=0}^n\b{np}{sp}\b{2sp}{sp}\b{2(n-s)p}{(n-s)p}
\\&\e \sum_{s=0}^n\b ns\b{2sp}{sp}\b{2(n-s)p}{(n-s)p}
 \e \sum_{s=0}^n\b ns\b{2s}s\b{2(n-s)}{n-s}
 =S_n\mod{p^{r+2}}.\endalign$$
Hence
$$\align &S_{np}-S_n\\&\e \sum_{t=1}^{p-1}\sum_{s=0}^{n-1}
\f{np}{sp+t}\b{(n-1)p+p-1}{sp+t-1}
\b{2sp+2t}{sp+t}\b{2(n-1-s)p+2(p-t)}{(n-1-s)p+p-t}\mod{p^{r+2}}.
\endalign$$
For $t\in\{1,2,\ldots,\f{p-1}2\}$ we have $\f p2<p-t<p$. By Lemma
2.6,
$$\b{2(n-1-s)p+2(p-t)}{(n-1-s)p+p-t}\e
(2(n-1-s)+1)\b{2(n-1-s)}{n-1-s}\b{2(p-t)}{p-t}\mod{p^2}.$$ By Lemma
2.2, $\b{2sp+2t}{sp+t}\e \b{2s}s\b{2t}t\mod p.$ Thus, applying Lemma
2.4 we see that
$$\align &\b{2sp+2t}{sp+t}
\b{2(n-1-s)p+2(p-t)}{(n-1-s)p+p-t}\\&\e \b{2s}s\b{2t}t
(2(n-1-s)+1)\b{2(n-1-s)}{n-1-s}\b{2(p-t)}{p-t}
\\&\e
-(2(n-1-s)+1)\b{2s}s\b{2(n-1-s)}{n-1-s}\f{2p}t\mod{p^2}.
\endalign$$
For $t\in\{\f{p+1}2,\ldots,p-1\}$ we have $1\le p-t<\f p2$. By Lemma
2.6,
$$\b{2sp+2t}{sp+t}\e (2s+1)\b{2s}s\b{2t}t\mod{p^2}.$$
By Lemma 2.2,
$$\b{2(n-1-s)p+2(p-t)}{(n-1-s)p+p-t}\e
\b{2(n-1-s)}{n-1-s}\b{2(p-t)}{p-t}\mod p.$$ Thus, applying Lemma 2.4
we get
$$\align &\b{2sp+2t}{sp+t}
\b{2(n-1-s)p+2(p-t)}{(n-1-s)p+p-t}\\&\e (2s+1)\b{2s}s\b{2t}t
\b{2(n-1-s)}{n-1-s}\b{2(p-t)}{p-t}
\\&\e
(2s+1)\b{2s}s\b{2(n-1-s)}{n-1-s}\f{2p}t\mod{p^2}.
\endalign$$
Hence
$$\align &S_{np}-S_n\\&\e
\sum_{t=1}^{(p-1)/2}\sum_{s=0}^{n-1}
\f{np}{sp+t}\b{(n-1)p+p-1}{sp+t-1}
\b{2sp+2t}{sp+t}\b{2(n-1-s)p+2(p-t)}{(n-1-s)p+p-t}
\\&\q+\sum_{t=(p+1)/2}^{p-1}\sum_{s=0}^{n-1}
\f{np}{sp+t}\b{(n-1)p+p-1}{sp+t-1}
\b{2sp+2t}{sp+t}\b{2(n-1-s)p+2(p-t)}{(n-1-s)p+p-t}
\\&\e -\sum_{t=1}^{(p-1)/2}\sum_{s=0}^{n-1}
\f{np}{sp+t}\b{(n-1)p+p-1}{sp+t-1}(2(n-1-s)+1)
\b{2s}s\b{2(n-1-s)}{n-1-s}
\f{2p}t
\\&\q+\sum_{t=(p+1)/2}^{p-1}\sum_{s=0}^{n-1}
\f{np}{sp+t}\b{(n-1)p+p-1}{sp+t-1}(2s+1)\b{2s}s\b{2(n-1-s)}{n-1-s}
\f{2p}t
\\&\e -\sum_{t=1}^{(p-1)/2}\sum_{s=0}^{n-1}
\f{np}t\b{n-1}s\b{p-1}{t-1}(2(n-1-s)+1)\b{2s}s\b{2(n-1-s)}{n-1-s}
\f{2p}t
\\&\q+\sum_{t=(p+1)/2}^{p-1}\sum_{s=0}^{n-1}
\f{np}t\b{n-1}s\b{p-1}{t-1}(2s+1)\b{2s}s\b{2(n-1-s)}{n-1-s} \f{2p}t
\\&\e
2np^2\sum_{s=0}^{n-1}(2(n-1-s)+1)
\b{n-1}s\b{2s}s\b{2(n-1-s)}{n-1-s}\sum_{t=1}^{(p-1)/2}\f{(-1)^t}{t^2}
\\&\q-2np^2\sum_{s=0}^{n-1}(2s+1)
\b{n-1}s\b{2s}s\b{2(n-1-s)}{n-1-s}\sum_{t=(p+1)/2}^{p-1}\f{(-1)^t}
{t^2}
\\&\e 2np^2\sum_{s=0}^{n-1}(2s+1)
\b{n-1}s\b{2s}s\b{2(n-1-s)}{n-1-s}
\\&\q\times\Big(\sum_{t=1}^{(p-1)/2}\f{(-1)^t}{t^2}-
\sum_{t=(p+1)/2}^{p-1}\f{(-1)^t} {t^2}\Big)\mod{p^{r+2}}.
\endalign$$
By Lemma 2.7,
$$\align &\sum_{s=0}^{n-1}(2s+1) \b{n-1}s\b{2s}s\b{2(n-1-s)}{n-1-s}
\\&=S_{n-1}+2\sum_{s=1}^{n-1}s\b{n-1}s\b{2s}s\b{2(n-1-s)}{n-1-s}
\\&=S_{n-1}+2(n-1)\sum_{s=1}^{n-1}\b{n-2}{s-1}\b{2s}s
\b{2(n-1-s)}{n-1-s}\\&=S_{n-1}+(n-1)S_{n-1}=nS_{n-1}. \endalign$$
Note that $B_{p-2}=0$ and $E_{2n}=2^{2n}E_{2n}\sls 12$. From Lemma
2.5 we see that
$$\sum_{k=1}^{(p-1)/2}\f{(-1)^k}{k^2}\e \f
12(-1)^{\f{p-1}2}E_{p-3}\Ls 12\e 2(-1)^{\f{p-1}2}E_{p-3}\mod p.$$
Thus,
$$\align &\sum_{k=1}^{(p-1)/2}\f{(-1)^k}{k^2}
 -\sum_{k=(p+1)/2}^{p-1}\f{(-1)^k}{k^2}
\\&=\sum_{k=1}^{(p-1)/2}\f{(-1)^k}{k^2}
 -\sum_{k=1}^{(p-1)/2}\f{(-1)^{p-k}}{(p-k)^2}
 \e 2\sum_{k=1}^{(p-1)/2}\f{(-1)^k}{k^2}
 \e 4(-1)^{\f{p-1}2}E_{p-3}\mod p.\endalign$$
 Now from the above we deduce that
 $S_{np}-S_n\e 2np^2\cdot nS_{n-1}\cdot 4(-1)^{\f{p-1}2}E_{p-3}
 \mod {p^{r+2}}.$
 This yields the result in this case.

\par Now assume $3\nmid n$. By Lemmas 2.9 and 2.11,
$$\b{3n}{3s}\e \b ns(1+9ns^2-9s)\mod{27},$$ and

$$\align\b {3n}{3s}\b{6s}{3s}\b{6(n-s)}{3(n-s)} &\e \b
ns\b{2s}s\b{2(n-s)}{n-s}(1+9n)(1+9ns^2-9s)\\
&\e\b ns\b{2s}s\b{2(n-s)}{n-s}(1+9n+9ns^2-9s) \mod{27}.\endalign$$
 Thus,
$$\sum_{s=0}^n\b{3n}{3s}\b{6s}{3s}\b{6(n-s)}{3(n-s)}
\e(1+9n)S_n+9\sum_{s=0}^n(ns^2-s)\b ns\b{2s}s\b{2(n-s)}{n-s}
\mod{27}$$ and so
$$\align &S_{3n}-S_n\\&\e\sum_{t=1}^{2}\sum_{s=0}^{n-1}
\f{3n}{3s+t}\b{3(n-1)+3-1}{3s+t-1}
\b{6s+2t}{3s+t}\b{6(n-1-s)+2(3-t)}{3(n-1-s)+3-t}
\\&+9nS_n+9\sum_{s=0}^n(ns^2-s)\b ns\b{2s}s\b{2(n-s)}{n-s}
\\&\e2n3^2nS_{n-1}(-5/4)+9nS_n
+9\sum_{s=0}^n(ns^2-s)\b ns\b{2s}s\b{2(n-s)}{n-s}
\\&\e9nS_n-9S_{n-1}+9\sum_{s=0}^n(ns^2-s)\b ns\b{2s}s\b{2(n-s)}{n-s}
\mod{27}.
\endalign$$
By (1.3), for $n\e 2\mod 3$ we have
$$S_n+S_{n-1}\e 4(3n(n+1)+1)S_n-32n^2S_{n-1}
=(n+1)^2S_{n+1}\e 0\mod 3,$$ for $n\e 1\mod 3$ we have
$$S_n-S_{n-1}\e n^2S_n-4(3n(n-1)+1)S_{n-1}=-32(n-1)^2
S_{n-2}\e 0\mod 3.$$ Thus, $$S_n\e \Ls n3 S_{n-1}\mod
3\qtq{for}n\not\e 0\mod 3.\tag 3.1$$
 Applying (3.1) and (2.2)
we have $$\align S_{3n}-S_n & \e
9(nS_n-S_{n-1})+9\sum_{s=0}^n(ns^2-s)\b
ns\b{2s}s\b{2(n-s)}{n-s}
\\&\e 9\sum_{s=0}^n(ns^2-s)\b
ns\b{2s}s\b{2(n-s)}{n-s}
\\&=9\sum_{s=0}^n\Big(ns(s-1)+\f{n-1}2(2s-n)+\f{n(n-1)}2\Big)
\b ns\b{2s}s\b{2(n-s)}{n-s}
\\&=9n\sum_{s=0}^ns(s-1)\b ns\b{2s}s\b{2(n-s)}{n-s}
+\f{9n(n-1)}2S_n\mod{27}.\endalign$$ If $s=3k+2$ for some
nonnegative integer $k$, using Lemma 2.2 we find that
$\b{2s}s=\b{3(2k+1)+1}{3k+2} \e \b{2k+1}k\b 12=0\mod 3$. Thus,
$3\mid s(s-1)\b{2s}s$ for any nonnegative integer $s$. Hence, from
the above we deduce that
$$S_{3n}-S_n\e \f{9n(n-1)}2S_n\e 9(n-n^2)S_n\e 9(n-1)S_n\mod{27}.$$
 This completes the proof.

 \pro{Corollary 3.1} Let $p>3$ be a prime. Then
 $$\align &S_p\e 4+8(-1)^{\f{p-1}2}p^2E_{p-3}\mod{p^3},
\\&S_{2p}\e 20+128(-1)^{\f{p-1}2}p^2E_{p-3}\mod{p^3},
\\&S_{3p}\e 112+1440(-1)^{\f{p-1}2}p^2E_{p-3}\mod{p^3}.\endalign$$
 \endpro
\pro{Remark 3.1} Let $p$ be an odd prime and $m,r\in\Bbb Z^+$. Then
$S_{m{p^r}} \equiv S_{m{p^{r-1}}}\pmod {p^{2r}}.$\endpro
 Proof.
 $S_{mp^r}=\sum\limits_{s=0}^{mp^{r-1}}\b{mp^r}{sp}
 \b{2sp}{sp}\b{2mp^r-2sp}{mp^r-sp}
+\sum\limits_{s=0}^{mp^{r-1}-1}\sum\limits_{t=1}^{p-1}
\b{mp^r}{sp+t}\b{2sp+2t}{sp+t}\b{2mp^r-2sp-2t}{mp^r-sp-t}.$
\\Applying Lemmas 2.12 and 2.14 we obtain
$$S_{mp^r}\e\sum_{s=0}^{mp^{r-1}}
\b{mp^{r-1}}s\b{2s}{s}\b{2mp^{r-1}-2s}{mp^{r-1}-s}=S_{mp^{r-1}}
\pmod{p^{2r}}.$$

\pro{Lemma 3.1} Let $m,n\in\Bbb Z^+$. Then
$$S_{mn+1}\e 4(mn+1)S_{mn}\mod{m^2n^2}.$$
\endpro
Proof. By (1.3),
$$(mn+1)^2S_{mn+1}=4(3mn(mn+1)+1)S_{mn}-32m^2n^2S_{mn-1}.$$
Thus,
$$(1+2mn)S_{mn+1}\e (mn+1)^2S_{mn+1}\e 4(3mn+1)S_{mn}
\mod{m^2n^2}$$ and so
$$S_{mn+1}\e \f {4(1+3mn)}{1+2mn}S_{mn}\e 4(1+3mn)(1-2mn)S_{mn}
\e 4(1+mn)S_{mn}\mod{m^2n^2}$$ as asserted.
 \pro{Theorem 3.2} Let
$p$ be an odd prime, and $m,r\in\Bbb Z^+$. Then
$$S_{mp^r+1}\e 4(mp^r+1)S_{mp^{r-1}}\mod{p^{2r}}.$$
\endpro

Proof. As $\t{ord}_p(m^2p^{2r})\ge 2r$, from Lemma 3.1 and Remark
3.1 we see that
$$S_{mp^r+1}\e 4(mp^r+1)S_{mp^r}\e 4(mp^r+1)S_{mp^{r-1}}
\mod {p^{2r}}.$$
This completes the proof.
 \pro{Theorem 3.3} Let $p$ be an odd prime and
$m,r\in\Bbb Z^+$. Then
$$S_{m{p^r}-1} \equiv (-1)^{\f {p-1}2}S_{m{p^{r-1}}-1}\pmod
{p^{r}}.$$\endpro
 Proof. It is clear that
 $$\aligned S_{mp^r-1}
 &=\sum_{s=0}^{mp^{r-1}-1}\sum_{t=0}^{p-1}
 \b{2sp+2t}{sp+t}\b{mp^r-1}{sp+t}\b{2(mp^r-1-sp-t)}{mp^r-1-sp-t}
\\&=\sum_{s=0}^{mp^{r-1}-1}\sum\limits\Sb {t=0}
\\t\neq(p-1)/2\endSb^{p-1}
 \b{2sp+2t}{sp+t}\b{mp^r-1}{sp+t}\b{2(mp^r-1-sp-t)}{mp^r-1-sp-t}
\\&+\sum_{s=0}^{mp^{r-1}-1}
\b{2sp+p-1}{sp+\f{p-1}2}\b{mp^r-1}{sp+\f{p-1}2}
\b{2(mp^r-1-sp-\f{p-1}2)}{mp^r-1-sp-\f{p-1}2} .\endaligned$$ Using
Lemma 2.8 we see that
$$\b{mp^r-1}{sp+t}\e\b{mp^{r-1}-1}{s}(-1)^t\pmod{p^r}.$$
For $t\neq\f{p-1}2$ applying Lemma 2.12 we obtain
$$\align&\b{2sp+2t}{sp+t}\b{2(mp^r-1-sp-t)}{mp^r-1-sp-t}
\\&=\b{2sp+2t}{sp+t}\b{2(mp^r-sp-t)}{mp^r-sp-t}
\f{(mp^r-sp-t)^2}{(2mp^r-1-2sp-2t)2(mp^r-sp-t)}
 \\&\e\b{2sp+2t}{sp+t}\b{2(mp^r-sp-t)}{mp^r-sp-t}
 \f{sp+t}{2(2sp+2t+1)}
 \e 0\pmod{p^r}.\endalign$$
For $t=\f{p-1}2$ using Lemma 2.13 we deduce that
$$\align S_{mp^r-1}&\e(-1)^{\f{p-1}2} \sum_{s=0}^{mp^{r-1}-1}
\b{2s}s\b{mp^{r-1}-1}s\b{2(mp^{r-1}-1-s)}{mp^{r-1}-1-s}
\\&=(-1)^{\f{p-1}2}S_{mp^{r-1}-1}\pmod{p^r}.\endalign$$
So the theorem is proved.

\pro{Theorem 3.4} Let $p$ be an odd prime and $n\in\Bbb Z^+$. Then
$$S_{np+1}
\e \cases (4+12n-9n^2)S_n\mod{p^3}&\t{if $p=3$,}
\\4(np+1)S_n+32n^2S_{n-1}(-1)^{\f{p-1}2}(E_{p-3}-1)p^2
\pmod {p^3}&\t{if $p>3$.}\endcases$$
\endpro
 Proof. By $(1.3)$,
 $(np+1)^2S_{np+1}=4(3np(np+1)+1)S_{np}-32n^2p^2S_{np-1}$.
  Thus, applying Theorems 3.1 and 3.3 we see that for $p>3$,
 $$\align (np+1)^2S_{np+1}&\e
 4(3n^2p^2+3np+1)(S_n+8n^2S_{n-1}(-1)^{\f{p-1}2}p^2E_{p-3})
 -32n^2p^2(-1)^{\f{p-1}2}S_{n-1}
 \\&\e 4(3n^2p^2+3np+1)S_n+32n^2S_{n-1}(-1)^{\f{p-1}2}(E_{p-3}-1)p^2
 \pmod {p^3},\endalign$$
and for $p=3$,
$$\aligned (3n+1)^2S_{3n+1}&=4(9n(3n+1)+1)S_{3n}-32n^2\cdot 9S_{3n-1}
\\&\e 4(9n+1)(1-9n(n-1))S_n-9n^2S_{n-1}
\\&\e 4(1-9n(n+1))S_n-9n^2S_{n-1}
\mod{27}.\endaligned\tag 3.2$$
 Since
 $$\f 1{(np+1)^2}=\f{(n^2p^2-np+1)^2}{\big((np)^3+1\big)^2}
 \e(n^2p^2-np+1)^2\e
 3n^2p^2-2np+1\mod {p^3},\tag 3.3$$
 from the above we deduce that for $p>3$,
 $$\aligned
 S_{np+1}&\e\f{4(3n^2p^2+3np+1)S_n+32n^2S_{n-1}
 (-1)^{\f{p-1}2}(E_{p-3}-1)p^2}
 {(np+1)^2}
\\&\e 4(S_n+3npS_n+n^2p^2(3S_n+8S_{n-1}(-1)^{\f{p-1}2}(E_{p-3}-1))
(3n^2p^2-2np+1)
\\&\e 4(np+1)S_n+32n^2S_{n-1}(-1)^{\f{p-1}2}(E_{p-3}-1)p^2\pmod
{p^3}.\endaligned$$
\par Now assume $p=3$. If $3\mid n$, from (3.2) and
(3.3) we deduce that
$$S_{3n+1}\e \f{4S_n}{(3n+1)^2}
\e 4(1-6n)S_n\e (4+12n-9n^2)S_n \mod{27}.$$ If $3\nmid n$, then
$S_{n-1}\e \sls n3S_n\mod 3$ by (3.1). Hence, from (3.2) and (3.3)
we deduce that
$$
\align S_{3n+1}&\e\f{4S_n-9n((n+1)S_n+nS_{n-1})}{(3n+1)^2} \e
\f{4S_n-9(n+1+\sls n3)S_n}{(3n+1)^2} \\&\e (4-9(2n+1))S_n(1-6n)\e
(12n-5)S_n\e (4+12n-9n^2)S_n \mod{27}.
\endalign$$
\par Summarizing the above proves the theorem.

 \pro{Corollary 3.2} Let $p>3$ be a
prime. Then
$$S_{p+1}\e 16+16p+32(-1)^{\f{p-1}2}(E_{p-3}-1)p^2\pmod {p^3}.$$
\endpro
 Proof. Taking $n=1$ in Theorem 3.4 we obtain the result.
\pro{Theorem 3.5} Let $p$ be a prime with $p\e 5,7\mod 8$. Then
$$S_{\f{p^2-1}2}\e 0\pmod{p^2}\qtq{and}
f_{\f{p^2-1}2}\e 0\pmod{p^2}.$$
\endpro
 Proof. For $\f{p-1}2<t< p$, from Lemma 2.2 we see that
 $p\mid \b{2sp+2t}{sp+t}$
  and $p\mid \b{\f{p-1}2p+\f{p-1}2}{sp+t}$. So
$$f_{\f{p^2-1}2}\e\sum_{s=0}^{\f{p-1}2}\sum_{t=0}^{\f{p-1}2}
 \b{\f{p-1}2p+\f{p-1}2}{sp+t}^3
 \pmod{p^2}$$ and

 $$S_{\f{p^2-1}2}\e\sum_{s=0}^{\f{p-1}2}\sum_{t=0}^{\f{p-1}2}
 \b{\f{p-1}2p+\f{p-1}2}{sp+t}\b{2sp+2t}{sp+t}
 \b{(p-1-2s)p+p-1-2t}{(\f{p-1}2-s)p+\f{p-1}2-t}
 \pmod{p^2}.$$
For $k\in\Bbb Z^+$ set $H_k=1+\f 12+\cdots+\f 1k$. For
$k\in\{1,\ldots,p-1\}$ we see that
$$\align\b{p-1}k&=\f{(p-1)(p-2)\cdots(p-k)}{k!}
\e\f{(-1)(-2)\cdots(-k)\big(1+\sum_{i=1}^k\f p{-i}\big)}{k!}
\\&=(-1)^k(1-pH_k)\pmod {p^2}\endalign$$ and so
$\f1{\b{p-1}{2s}}\e 1+pH_{2s}\pmod {p^2}$ for
$s=1,2,\ldots,\f{p-1}2$.
 It is well known that $H_{\f{p-1}2}\e-\f{2^p-2}p \pmod p$. Thus,
$$\b{p-1}{\f{p-1}2}\e(-1)^{\f{p-1}2}(1-pH_{\f{p-1}2})
\e(-1)^{\f{p-1}2}(2^p-1)\pmod {p^2}. $$
 Applying (2.4) and Lemma 2.8
we have
 $$\aligned&\b{2sp+2t}{sp+t}\b{p^2-1-2sp-2t}{\f{p^2-1}2-sp-t}
\\&=\f{\b{p^2-1}{\f{p^2-1}2}\b{\f{p-1}2p+\f{p-1}2}{sp+t}^2}
{\b{p^2-1}{2sp+2t}}
=\f{\b{p-1}{\f{p-1}2}\b{\f{p-1}2p+\f{p-1}2}{sp+t}^2}{\b{p-1}{2s}}
 \f{\prod\limits\Sb i=1\\p\nmid i\endSb ^{(p^2-1)/2}\f{p^2-i}i}
 {\prod\limits\Sb i=1\\p\nmid i\endSb ^{2sp+2t}\f{p^2-i}i}
\\&\e(-1)^{\f{p-1}2}(2^p-1)(1+pH_{2s})\b{\f{p-1}2p+\f{p-1}2}{sp+t}^2
\\&\e\Big((-1)^{\f{p-1}2}(2^p-1)+(-1)^{\f{p-1}2}pH_{2s})\Big)
\b{\f{p-1}2p+\f{p-1}2}{sp+t}^2\pmod {p^2} \endaligned$$ and so
$$ S_{\f{p^2-1}2}
 \e(-1)^{\f{p-1}2}(2^p-1)f_{\f{p^2-1}2}
 +(-1)^{\f{p-1}2}p\sum_{s=0}^{(p-1)/2}H_{2s}\sum_{t=0}^{(p-1)/2}
 \b{\f{p-1}2p+\f{p-1}2}{sp+t}^3
 \pmod{p^2}.$$
\par Now we assert that
 $$\sum\limits_{t=0}^{(p-1)/2}\b{\f{p-1}2p+\f{p-1}2}{sp+t}^3\e 0
  \pmod {p^2} \qtq {for} s=0,1,2,\ldots.\tag 3.4$$
  We prove the result by induction on s.
  From [S1] we know that the result is true for $s=0$.
 Suppose
 $$\align&\sum\limits_{t=0}^{(p-1)/2}\b{\f{p-1}2p+\f{p-1}2}{sp+t}^3
 \\&\e\b {\f{p-1}2}s^3\sum\limits_{t=0}^{(p-1)/2}
 \Big[\f{p+1}2\b {\f{p-1}2}t-\big(\f{p-1}2+s\big)\b {\f{p-1}2-p}t-
 s\b {\f{p-1}2-p}{t+p}\Big]^3
 \e 0
  \pmod {p^2}.\endalign$$
 By Lemma 2.15,
 $$\aligned&\sum\limits_{t=0}^{(p-1)/2}
 \b{\f{p-1}2p+\f{p-1}2}{(s+1)p+t}^3
 \\&\e\b {\f{p-1}2}{s+1}^3\sum\limits_{t=0}^{(p-1)/2}
 \Big[\f{p+1}2\b {\f{p-1}2}t-\big(\f{p-1}2+s+1\big)\b {\f{p-1}2-p}t-
 (s+1)\b {\f{p-1}2-p}{t+p}\Big]^3
 \\&=\b {\f{p-1}2}{s+1}^3\sum\limits_{t=0}^{(p-1)/2}
 \Big[\f{p+1}2\b {\f{p-1}2}t-\big(\f{p-1}2+s\big)\b {\f{p-1}2-p}t
 \\&\q-s\b {\f{p-1}2-p}{t+p}
 -\Big(\b {\f{p-1}2-p}t+\b {\f{p-1}2-p}{t+p}\Big)\Big]^3
\pmod {p^2}.\endaligned $$
 Hence $\sum\limits_{t=0}^{(p-1)/2}
 \b{\f{p-1}2p+\f{p-1}2}{(s+1)p+t}^3\e 0\mod{p^2}$ for $s\ge \f{p-1}2$.
 For $s<\f{p-1}2$,  by the inductive hypothesis we have
 $$\sum\limits_{t=0}^{(p-1)/2}
 \Big[\f{p+1}2\b {\f{p-1}2}t-(\f{p-1}2+s)\b {\f{p-1}2-p}t-
 s\b {\f{p-1}2-p}{t+p}\Big]^3
 \e 0\pmod {p^2}.$$
 As $t\in\{0,1,\ldots,\f{p-1}2\}$, we have
 $\b {\f{p-1}2}t\e\b {\f{p-1}2-p}t\e\b {-\f{1}2}t\pmod p.$
 From Lemma 2.16 we  have
 $\b{\f{p-1}2-p}t+\b{\f{p-1}2-p}{t+p}
 =(-1)^t\Big(\b{\f{p-1}2+t}t-\b{p+\f{p-1}2+t}{t+p}\Big)
 \e 0\pmod p,$
  and so
 $$\align&\sum\limits_{t=0}^{(p-1)/2}
 \b{\f{p-1}2p+\f{p-1}2}{(s+1)p+t}^3
\\&\e\b {\f{p-1}2}{s+1}^3\{\sum\limits_{t=0}^{(p-1)/2}
 \Big[\f{p+1}2\b {\f{p-1}2}t-\big(\f{p-1}2+s\big)\b {\f{p-1}2-p}t-
 s\b {\f{p-1}2-p}{t+p}\Big]^3
\\&\q+3\sum\limits_{t=0}^{(p-1)/2}
 \Big[\f{p+1}2\b {\f{p-1}2}t-\big(\f{p-1}2+s\big)\b {\f{p-1}2-p}t-
 s\b{\f{p-1}2-p}{t+p}\Big]
 \\&\q\q\q\q\q\times\Big[\b{\f{p-1}2-p}t+\b{\f{p-1}2-p}{t+p}\Big]^2
\\&\q -3\sum\limits_{t=0}^{(p-1)/2}
 \Big[\f{p+1}2\b {\f{p-1}2}t-\big(\f{p-1}2+s\big)\b {\f{p-1}2-p}t-
 s\b{\f{p-1}2-p}{t+p}\Big]^2
 \\&\q\q\q\q\q\times\Big[\b{\f{p-1}2-p}t+\b{\f{p-1}2-p}{t+p}\Big]
 \\&\q-\sum\limits_{t=0}^{(p-1)/2}
 \Big[\b{\f{p-1}2-p}t+\b{\f{p-1}2-p}{t+p}\Big]^3\}
\\&\e -3\b{\f{p-1}2}{s+1}^3\sum\limits_{t=0}^{(p-1)/2}
 \b{-\f{1}2}t^2\Big[\b{\f{p-1}2-p}t+\b{\f{p-1}2-p}{t+p}\Big]
 \\&=-3\b{\f{p-1}2}{s+1}^3\sum\limits_{t=0}^{(p-1)/2}
 (-1)^t\b{-\f{1}2}t^2
 \Big[\b{\f{p-1}2+t}t-\b{p+\f{p-1}2+t}{t+p}\Big]\\&= 0\pmod {p^2}.
 \endalign $$
\par Summarizing the above proves the theorem.

\section*{4. $\{S_m\}$ is log-convex}
\par A sequence $\{a_m\}$ $(m\ge 0)$ is called log-convex if $a_m\ge
0$ and $a_{m-1}a_{m+1}\ge a_m^2$ for $m=1,2,3,\ldots$. In this
section we show that $\{S_m\}$ and $\{P_m\}$ are log-convex
sequences.
 \pro{Theorem 4.1}
For $m=2,3,4,\ldots$ we have
$$S_m^2<S_{m+1}S_{m-1}<\big(1+\f 1{m(m-1)}\big)S_m^2.$$
\endpro
Proof. We first prove $S_m^2<S_{m+1}S_{m-1}$ for $m\ge 2$. Since
$S_1=4$, $S_2=20$, $S_3=112$, $S_4=676$ and $S_5=4304$ we see that
$S_m^2<S_{m+1}S_{m-1}$ for $m=2,3,4$. From now on we assume $m\ge
5$. Suppose that $S_{m-1}^2<S_{m-2}S_m$.
  By (1.3), Lemma 2.7 and the fact that $\f{3m-4}{m-1}\ge \f {11}4$
  we have
  $$\align &S_{m+1}S_{m-1}-S_m^2
 \\&=\f{4(3m^2+3m+1)}{(m+1)^2}S_mS_{m-1}
 -\f{32m^2}{(m+1)^2}S_{m-1}^2-S_m^2
 \\&>\Big(\f{4(3m^2+3m+1)}{(m+1)^2}S_{m-1}
 -\f{32m^2}{(m+1)^2}S_{m-2}-S_m\Big)S_m
\\&= \Big(4(3m^2+m-1)S_{m-1}
+32(1-2m^2)S_{m-2}\Big)\f{S_m}{m^2(m+1)^2}
\\&=\Big(4(3m^2+m-1)2\sum_{k=1}^{m-1}
\b{m-2}{k-1}\b{2k}{k}\b{2m-2-2k}{m-1-k}
\\&\q+32(1-2m^2)\sum_{k=0}^{m-2}
\b{m-2}k\b{2k}{k}\b{2m-4-2k}{m-2-k}\Big)\f{S_m}{m^2(m+1)^2}
\\&=\Big(\sum_{k=1}^{m-3}\b{m-2}k\b{2k}{k}\b{2m-4-2k}{m-2-k}
\Big(\f{16(3m^2+m-1)(2k+1)}{k+1}+32(1-2m^2)\Big)
\\&\q+\Big(16(3m^2+m-1)\f{3m-4}{m-1}+64-128m^2\Big)
\b{2m-4}{m-2}\Big)\f{S_m}{m^2(m+1)^2}
\\&>\Big(8(m^2+3m+1)\sum_{k=1}^{m-3}
\b{m-2}k\b{2k}{k}\b{2m-4-2k}{m-2-k}
\\&\q+4(m^2+11m+5)\b{2m-4}{m-2}\Big)
\f{S_m}{m^2(m+1)^2}\\&>0.\endalign$$ Thus the inequality
$S_m^2<S_{m+1}S_{m-1}$ is proved by induction.
 \par Next we prove
the remaining inequality.
 It is easily seen that
  $\big(1+\f 1{m(m-1)}\big)S_m^2-S_{m+1}S_{m-1}>0$ for
  $m=2,3,\ldots,13$. Now suppose $m\ge 14$ and
  $\big(1+\f 1{(m-1)(m-2)}\big)S_{m-1}^2>S_{m}S_{m-2}$.
By (1.3), Lemma 2.7 and the inductive hypothesis we have
  $$\align &\big(1+\f 1{m(m-1)}\big)S_m^2-S_{m+1}S_{m-1}
  \\&=\big(1+\f
  1{m(m-1)}\big)S_m^2-\f{4(3m^2+3m+1)}{(m+1)^2}S_mS_{m-1}
  +\f{32m^2}{(m+1)^2}S_{m-1}^2
\\&>\Big(\f{m^2-m+1}{m(m-1)}S_m-\f{4(3m^2+3m+1)}{(m+1)^2}S_{m-1}
+\f{32m^2(m-1)(m-2)}{(m+1)^2(m^2-3m+3)}S_{m-2}\Big)S_m
\\&=\Big((20m^5-60m^4+52m^3+28m^2-36m+12)S_{m-1}
\\&\q+(-128m^5+320m^4-256m^3-32m^2+192m-96)S_{m-2}\Big)
\\&\q\q\times\f{S_m}{(m+1)^2(m^2-3m+3)m^3(m-1)}
\\&=16\sum_{k=0}^{m-2}\b{m-2}k\b{2k}k\b{2m-4-2k}{m-2-k}
F(m,k)\f{S_m}{(m+1)^2(m^2-3m+3)m^3(m-1)},
\endalign $$
where $$\align F(m,k)&=(5m^5-15m^4+13m^3+7m^2-9m+3)
\f{2k+1}{k+1}\\&\q-8m^5+20m^4-16m^3-2m^2+12m-6.\endalign$$
 For $m\ge14$  we see that $3<\f{(2m-7)(2m-5)}{(m-3)(m-2)}<4$,
 $5m^5-15m^4+13m^3+7m^2-9m+3>0$,
 $-8m^5+20m^4-16m^3-2m^2+12m-6<0$,
 $6m^7-75m^6+223m^5-283m^4-61m^3+427m^2-87m-42>0,$
 and  $F(m,k+1)>F(m,k)$ for $k=0,1,\ldots,m-2$. Thus
$F(m,m-3)+F(m,1)>F(m,5)+F(m,1)>0.$ Since
$$\align F(m,k)\ge F(m,2)&= \f 53(5m^5-15m^4+13m^3+7m^2-9m+3)
\\&\q-8m^5+20m^4-16m^3-2m^2+12m-6>0,\endalign$$
from the above we derive that
$$\align &\q\big(1+\f 1{m(m-1)}\big)S_m^2-S_{m+1}S_{m-1}
\\&>16\Big\{\sum_{k=0}^{2}\b{m-2}k\b{2k}k\b{2m-4-2k}{m-2-k}F(m,k)
\\&\q+\sum_{k=m-4}^{m-2}\b{m-2}k\b{2k}k\b{2m-4-2k}{m-2-k}F(m,k)\Big\}
\f{S_m}{(m+1)^2(m^2-3m+3)m^3(m-1)}
\\&=\f{16S_m}{(m+1)^2(m^2-3m+3)m^3(m-1)}
\Big\{\b{2m-4}{m-2}\big(F(m,m-2)+F(m,0)\big)
\\&\q+2(m-2)\b{2m-6}{m-3}\big(F(m,1)+F(m,m-3)\big)
\\&\q+3(m-2)(m-3)\b{2m-8}{m-4}\big(F(m,m-4)+F(m,2)\big)\Big\}
\\&>\f{16S_m}{(m+1)^2(m^2-3m+3)m^3(m-1)}
\\&\q\times\Big(\big(3(m^2-5m+6)+\f{4(2m-7)(2m-5)}{(m-3)(m-2)}\big)F(m,m-4)
\\&\q+\f{4(2m-7)(2m-5)}{(m-3)(m-2)}F(m,0)\Big)\b{2m-8}{m-4}
\\&>\f{16S_m\b{2m-8}{m-4}}{(m+1)^2(m^2-3m+3)m^3(m-1)}
\\&\q\times\Big(3(m^2-5m+6)F(m,m-4)
+\f{4(2m-7)(2m-5)}{(m-3)(m-2)}F(m,0)\Big)
\\&=\f{S_m\b{2m-8}{m-4}}{(m+1)^2(m^2-3m+3)m^3(m-1)}
\Big(\big(6(m-2)(2m-7)+\f{8(2m-7)(2m-5)}{(m-3)(m-2)}\big)
\\&\q\q\q\q\q\q\q\q\q\q\q\q\q\q\q\q\q\times(40m^5-120m^4+104m^3+56m^2-72m+24)
\\&\q+\big(3(m-2)(m-3)+\f{4(2m-7)(2m-5)}{(m-3)(m-2)}\big)
\\&\q\q\times(-128m^5+320m^4-256m^3-32m^2+192m-96)\Big)
\\&>\f{S_m\b{2m-8}{m-4}}{(m+1)^2(m^2-3m+3)m^3(m-1)}
\\&\q\times\Big(\big(6(m-2)(2m-7)+24\big)
(40m^5-120m^4+104m^3+56m^2-72m+24)
\\&\q+\big(3(m-2)(m-3)+16\big)(-128m^5+320m^4-256m^3-32m^2+192m-96)\Big)
\\&=\f{S_m\b{2m-8}{m-4}}{(m+1)^2(m^2-3m+3)m^3(m-1)}
\\&\q\times(96m^7-1200m^6+3568m^5-4528m^4-976m^3+6832m^2-1392m-672)
\\&>0.\endalign$$ Hence the inequality is proved by induction.

\pro{Corollary 4.1} Both $\{S_m\}$ and $\{P_m\}$ are log-convex.
\endpro
Proof. By (1.6), $S_m>0$. Since $S_0=1$, $S_1=4$ and $S_2=20$, we
see that $S_1^2<S_0S_2$. Now applying Theorem 4.1 we see that
$\{S_m\}$ is log-convex. Since $P_m=2^mS_m$ we see that
$$P_{m-1}P_{m+1}-P_m^2=2^{m-1}S_{m-1}\cdot
2^{m+1}S_{m+1}-2^{2m}S_m^2=2^{2m}(S_{m-1}S_{m+1}-S_m^2)\ge 0.$$ Thus
 $\{P_m\}$ is also log-convex.

\end{document}